

\input amstex
\input amsppt.sty
\magnification=\magstep1
\hsize = 6.5 truein
\vsize = 9 truein

\NoRunningHeads
\NoBlackBoxes
\TagsAsMath

\def\label#1{\par%
        \hangafter 1%
        \hangindent .75 in%
        \noindent%
        \hbox to .75 in{#1\hfill}%
        \ignorespaces%
        }

\newskip\sectionskipamount
\sectionskipamount = 24pt plus 8pt minus 8pt
\def\sectionskip{\vskip\sectionskipamount}
\define\sectionbreak{%
        \par \ifdim\lastskip<\sectionskipamount
        \removelastskip    \penalty-2000 \sectionskip \fi}
\define\section#1{%
        \sectionbreak 
        \subheading{#1}%
        \bigskip
        }

\redefine\qed{{\unskip\nobreak\hfil\penalty50\hskip2em\vadjust{}\nobreak\hfil
    $\square$\parfillskip=0pt\finalhyphendemerits=0\par}}

\define\op#1{\operatorname{\fam=0\tenrm{#1}}} 

        \let \< = \langle \let \> = \rangle \define \a {\alpha}
\define \e {\varepsilon} \define \G {\Gamma} \define \h {\hslash}
\define \th {\theta} \define \p {\partial}

\topmatter
\title The Classical Limit of Dynamics for Spaces Quantized by an
  Action of $\Bbb R^d$
\endtitle
\author Marc A. Rieffel
\endauthor
\thanks The research reported here
   was supported in part by National Science Foundation grant
DMS--9303386.
\endthanks
\abstract
  We have previously shown how to construct a deformation quantization
of any locally compact space on which a vector group acts. Within this
framework we show here that, for a natural class of Hamiltonians, the
quantum evolutions will have the classical evolution as their
classical limit.
\endabstract
\address\hskip-\parindent
        Department of Mathematics \newline University of California
\newline Berkeley, California 94720-3840
\endaddress
\email rieffel\@math.berkeley.edu
\endemail
\subjclass   Primary 46L60; Secondary 46l55 81S30
\endsubjclass
\endtopmatter

\baselineskip=14 pt
\document

\section{Introduction}

Let $M$ be a locally compact space, and let $\a$ be an action of
$V=\Bbb R^d$ on $M$. Let $A=C_\infty(M)$, the C${}^*$-algebra of
complex-valued continuous functions on $M$ which vanish at infinity,
and let $\a$ denote also the corresponding action of $V$ on $A$.  Let
$J$ be a skew-symmetric operator on $V$. Then $J$ determines a
``Poisson bracket'' on $A$, and in \cite{Rf} we have shown how to
construct a strict deformation quantization of $A$ into a
one-parameter family, $A_{\hslash}$, of non-commutative
C${}^*$-algebras, in ``the direction of this Poisson bracket''. The
purpose of the present paper is to show that within this framework,
the quantum evolution of the system which is determined by any
Hamiltonian from a natural class, converges as $\hslash\to 0$ to the
classical evolution for that Hamiltonian, as one would expect.

Our aim here is not at all to obtain the strongest possible results
--- a more elaborate and lengthy analysis could deal with a far wider
class of Hamiltonians than the ones we consider here. Rather our aim
is to show how naturally this matter fits within the framework of
\cite{Rf}. The main argument, given in section 2, is relatively simple
and brief (though heavily dependent on the results in \cite{Rf}).

The rest of this paper, contained in sections 3 and 4, is concerned
simply with showing that the Hamiltonians which we consider do have
classical flows (evolutions) which exist for all time, and that these
classical flows have the smoothness properties needed for our analysis
in section 2. This is a matter of some independent interest, as
indicated in \cite{B}, and our situation permits us to remove the
restriction to locally free actions which is required in some of the
relevant places in \cite{B}. (I thank Alan Weinstein for a suggestion
which simplified my proof of the existence of the classical flow.)
This gives a partial answer to question 2.4.29 in \cite{B}.  In the
appendix we sketch how to obtain the existence of our classical flows
as a consequence of a powerful theorem of D. Robinson \cite{Rs1}. (I
thank George Elliott and Ola Bratteli for comments which led me to
look at this paper of Robinson.) It is natural to carry this out in
the more general context of an arbitrary Lie group acting on a space.
This gives a more substantial partial answer to question 2.4.29 of
\cite{B}. But I have not seen how to use this approach to conveniently
give the smoothness properties which we need for the proof of our main
theorem (see question 1 of \cite{Rs1}), and so I have found it best to
include the much more elementary approach given in section 3, since it
develops most of the tools needed for section 4.

There is already an enormous literature concerned with the classical
limit of quantum evolutions, mostly on $\Bbb R^{2n}$, and we will not
try to review it here. Many references can be found by chasing back
the references given in \cite{E, Rr, W}.

The construction of strict deformation quantizations developed in
\cite{Rf} works equally well for non-commutative C${}^*$-algebras, and so
one can ask whether the results of the present paper extend to that
case. The difficulty is that usually the Poisson bracket applied to a
Hamiltonian does not give a derivation of the non-commutative algebra,
and so one cannot expect it to generate a group of automorphisms
analogous to the classical flow. In other words, I don't know how to
even pose the question we consider here, for the more general
situation. (In the very special case where all is sufficiently related
to the center of the algebra one will obtain a derivation, and
presumably the results of the present paper can be extended to that
case; but at present it is not clear to me that this is of any
particular interest and so I have not pursued it here.)

\section{1. The classical flows}

The purpose of this section is to describe the classical vector
fields, and corresponding classical flows, which we will consider, and
to state those properties of these classical flows which we will need
later. We will defer the proofs of these properties until after our
discussion of the classical limits of the quantum flows in the next
section (where we will also relate our classical vector fields to
Hamiltonians).

Let $M$, $V$, $A$ and $\a$ be as in the introduction. It is the action
$\a$ which gives $M$ its ``differential'' structure.  (It would
certainly be of interest to consider actions of more general Lie
groups than $V$, but I don't know how to construct deformation
quantizations in that generality.)  Let $C_b(M)$ denote the algebra of
bounded continuous functions on $M$, the multiplier algebra of $A$.
The evident action $\a$ of $V$ on $C_b(M)$ is not in general strongly
continuous. Let $B$ (or $B(M)$, or $B(M,\a)$) denote the subspace of
elements of $C_b(M)$ on which $\a$ is norm-continuous, so that $B$ is
the largest C${}^*$-subalgebra of $C_b(M)$ on which $\a$ is strongly
continuous. Note that $B$ is a unital C${}^*$-algebra containing $A$
as an essential ideal.  (We could now view $B$ as the algebra of
continuous functions on its maximal ideal space, which is compact and
on which $\a$ gives an action, but this is not technically
advantageous at this point.)

Let $A^\infty$ and $B^\infty$ denote the dense $*$-subalgebras of $A$
and $B$ consisting of smooth (i.e. infinitely differentiable) vectors
\cite{B} for $\a$. As in chapter 9 of \cite{Rf} we will distinguish
between $V$ and its Lie algebra, denoting the latter by $L$.  Thus for
each $X\in L$ we have a corresponding derivation, $\a_X$, on
$A^\infty$ and $B^\infty$, given by the infinitesimal generator of the
one-parameter group of operators corresponding to $X$. We will
heuristically think of $\a_X$ as a smooth tangent vector field on $M$,
and think of the tangent space at each point $m$ of $M$ as
corresponding to $L$, by means of the $\a_X$'s followed by a point
evaluation at $m$. Then we will think of a continuous real
vector-field on $M$ as being just a continuous $L$-valued function on
$M$.

Fix for the rest of this paper an arbitrary positive-definite inner
product on $L$. Heuristically this makes $M$ into a Riemannian
``manifold''.  Let $C_b(M,L)$ denote the Banach space of continuous
bounded $L$-valued functions on $M$, equipped with the supremum norm
using the inner product on $L$. We have the evident action $\a$ of $V$
on $C_b(M,L)$. Much as above, denote by $B(M,L)$ the largest subspace
of $C_b(M,L)$ on which $\a$ is strongly continuous, and by
$B^\infty(M,L)$ the subspace of smooth vectors for $\a$. We think of
$B^\infty(M,L)$ as the space of smooth bounded vector fields on $M$.

Let $\Phi\in B^\infty(M,L)$. For any $f\in B^\infty$ and any $m\in M$
the function $x\mapsto f(\a^{-1}_x(m))$ is smooth on $V$, and in
particular it will have a finite total derivative, $(Df)_m$, at $x=0$.
This is a linear functional on $L$. Thus we can define a map,
$\delta_\Phi$, from $B^\infty$ to itself by $$ (\delta_\Phi f)(m) =
(Df)_m (\Phi(m)) \ .  \tag 1.1 $$ Then $\delta_\Phi$ is a
$*$-derivation of $B^\infty$, in accordance with our heuristic view
that $\Phi$ is a smooth vector field.  To see this, note that the
above is just a coordinate-free way of saying the following. Let $\{
E_j\}$ be a basis for $L$, and let $\{\Phi_j\}$ denote the
corresponding components of $\Phi$.  Note that $\Phi_j\in B^\infty$
for each $j$. Let $\p_j=\a_{E_j}$, a $*$-derivation of $B^\infty$.
Then $$
\delta_\Phi f = \Sigma \,\Phi_j (\p_j f)    \tag 1.2
$$ for $f\in B^\infty$. It is now clear that $\delta_\Phi$ is a
$*$-derivation of $B^\infty$, and that it carries $A^\infty$ into
itself.  (If $A$ were non-commutative, we could not expect
$\delta_\Phi$ to be a derivation unless each $\Phi_j$ were in the
center of $B$.)

The main fact which we need is that each $\Phi\in B^\infty(M,L)$
determines a flow on $M$ which exists for all time, and which carries
$B$ into itself. We formulate this as:

\proclaim{1.3 THEOREM} Let $M$, $\a$ and $B$ be as above, and
let $\Phi\in B^\infty(M,L)$.  Then $\delta_\Phi$ is a pregenerator,
that is, there is a (unique) strongly continuous one-parameter group,
$\beta$, of automorphisms of $B$ whose generator is the closure of
$\delta_\Phi$. Furthermore, $\beta$ carries $A$ into itself, and so
$\beta$ comes from a flow on $M$ (which we will also denote by
$\beta$).
\endproclaim

We will defer the proof of this theorem to section 3.

We also need control over the higher derivatives associated with the
flow $\beta$. For each $f\in B^\infty$ we have the higher total
derivatives, $D^kf$, of $f$. Thus each $D^kf$ is a function on $M$
into the (symmetric) $k$-linear functionals from $L$ to the complex
numbers, defined by $$ (D^kf)_m (X_1,\dots ,X_k) = (\a_{X_1}\dots
\a_{X_k}f)(m) $$ for each $m\in M$. Each $D^kf$ is smooth and bounded,
because of the definition of $B^\infty$. Thus, by using the inner
product on $L$ to define the norm of $k$-linear functionals on $L$, we
can define semi-norms $\| \ \ \|_{(k)}$ on $B^\infty$ by $$
\|f\|_{(k)} = \|D^kf\|_\infty =\sup
\{\| (D^kf)_m\|: m\in M\} \ .
$$ We will need:

\proclaim{1.4 THEOREM} With notation as above, let $\beta$ be the
action on $B^\infty$ for $\Phi\in B^\infty(M,L)$ as in Theorem 1.3.
Then the action $\beta$ is strongly continuous for each of the
semi-norms $\| \ \ \|_{(k)}$ on $B^\infty$. Furthermore, for any $f\in
B^\infty$ the function $t\mapsto\beta_t f$ is smooth for these
semi-norms, and its first derivative is $(\delta_\Phi f)\circ\beta_t$.
\endproclaim

We remark that $\beta$ will not usually be uniformly bounded for the
above semi-norms.

We defer the proof of this theorem to section 4.

\bigskip
\section{2. The classical limit}

As above, we let $\a$ be an action of $V$ on a locally compact space
$M$. Thus we have the algebras $A$ and $B$, and their smooth versions
$A^\infty$ and $B^\infty$. We let $J$ be a skew-symmetric operator on
$L$, so that $J$ determines a Poisson bracket, $\{ \ , \ \}$, on
$A^\infty$ and $B^\infty$. It is defined, in terms of a basis $\{
E_j\}$ for $L$, by $$
\{ f,g\} = \sum J_{jk}\a_{E_j}(f)\a_{E_k}(g) \ .
$$ For each ``Planck's constant'' $\h$ we let $A_\h$ and $B_\h$ denote
the corresponding deformed C${}^*$-algebras, as constructed in
\cite{Rf}.  Thus $A_{\h}$ has $A^\infty$ as dense subspace, with
product given there by $$ f\times_{\h} g =\iint \a_{\h Ju}(f)\a_v (g)
e(u\cdot v) \ du \ dv $$ (an oscillatory integral, with $e(t)=e^{2\pi
it}$), and with corresponding C${}^*$-norm. The involution is still
complex-conjugation.  We define $B_\hslash$ similarly. Then $A_{\h}$
will be an essential ideal in $B_\h$ by proposition 5.9 of \cite{Rf}.
Furthermore, $\a$ gives an action of $V$ on $A_\h$ and $B_\h$ by
proposition 5.11 of \cite{Rf}, and the corresponding subspaces of
smooth vectors will be exactly $A^\infty$ and $B^\infty$ as vector
spaces, by theorem 7.1 of \cite{Rf}.

The Hamiltonians which we will consider consist of the real-valued
functions in $B^\infty$. So fix such a real-valued $H\in B^\infty$.
The mapping $f\mapsto\{ H,f\}$ is a derivation of $A^\infty$ and
$B^\infty$. If for a basis $\{E_j\}$ for $L$ we set $$
\Phi(m) = \sum J_{jk} (\a_{E_j}(H))(m)E_k \ ,
$$ we see that $\Phi$ is a function from $M$ to $L$ of the kind
considered in the first section. In particular its coefficient
functions are in $B^\infty$ (and real) so that $\Phi\in
B^\infty(M,L)$. Furthermore, it is clear that $\{ H,f\}=\delta_\Phi
(f)$ for each $f\in B^\infty$. Thus $\Phi$ is the ``Hamiltonian vector
field'' for $H$.

According to Theorem 1.3, $\delta_\Phi$ determines a flow, $\beta$, on
$M$, with corresponding strongly continuous one-parameter action on
$A$ and $B$. This is the ``Hamiltonian flow'' for $H$.

For each $\h$ we let $[ \ , \ ]_\h$ denote the ordinary commutator for
the corresponding product in $B^\infty$, so that $$ [f,g]_\h =
f\times_\h g-g\times_\h f \ .  $$ Set $H^\h=(-\pi/ \h)H$, viewed as a
self-adjoint element of $B_\h$.  (The $-\pi$ comes from our
conventions in \cite{Rf} for the definition of $\times_\h$, given
above.) Then the map $f\mapsto [iH^\h,f]_\h$ is a $*$-derivation of
$A_\h$ and $B_\h$ which is bounded (but with norms going to $+\infty$
as $\h\to 0$).  Thus $H^\h$ determines a one-parameter group,
$\beta^\h$, of $*$-automorphisms of $A_\h$ and $B_\h$, the
corresponding quantum flow. This flow consists of inner automorphisms
of $B_\h$.  For let $u^\h_t=\exp_\h(it H^\h)$ for each $t\in\Bbb R$,
where $\exp_\h$ denotes the exponential defined by the usual power
series, but using the product $\times_\h$ in $B_\h$. Then $u_t^\h$ is
a unitary element of $B_\h$. By the usual calculations we will have $$
\beta^\h_t(f) = u^\h_t \times_\h f \times_\h u_{-t}^\h
$$ for all $f\in B_\h$. Again it is clear that $\beta^\h$ carries the
ideal $A_\h$ into itself. But notice that $\beta^\h$ is not only
strongly continuous, but actually norm (i.e. uniformly) continuous
(since $H^\h$ is bounded).

The main theorem of this paper is:

\proclaim{2.1 THEOREM} With notation as above, for any $f\in B^\infty$
we have $$
\| \beta^\h_t f-\beta_t f\|_\h \to 0 \ \ \ {\text{ as }} \ \ \
\h\to 0
$$ for each $t\in\Bbb R$, with the convergence being uniform in $t$
over any finite interval.
\endproclaim

It is in this sense that, within our framework, the quantum flow has
the classical flow as its classical limit.

We remark that in the proof we will see how to obtain specific
estimates for the convergence.

\bigskip
{\smc Proof}. Let $I$ denote the interval $[-1,1]$. We only need
consider $\h$'s in $I$. We will denote by $B^{(k)}$ the space
$B^\infty$ equipped with the norm $\| \ \ \|_k$ which is the sum of
$\| \ \ \|_\infty$ with the semi-norms $\| \ \ \|_{(j)}$ (defined near
the end of section 1) for $j\leq k$. This norm is equivalent to the
norm used in \cite{Rf}, defined on page 1 of \cite{Rf}.  We choose $k$
large enough that we can apply the little argument near the beginning
of the proof of theorem 9.3 of \cite{Rf} which shows that there is a
constant, $c$, independent of $f\in B^\infty$, such that $$
\| f\|_\h \leq c\|f\|_k   \tag 2.2
$$ for all $\h\in I$. Fix $f\in B^\infty$. From Theorem 1.4 we know
that $t\mapsto\beta_t f$ can be viewed as a smooth function on $V$
with values in $B^{(k)}$, whose first derivative is $\{ H,\beta_t
f\}$.  From 2.2 it follows that $t\mapsto\beta_t^f$ is smooth as a
function with values in $B_\h$, for each $\h\in I$, with the same
first derivative.

Fix $\h\in I$. Then the smooth function $t\mapsto u^\h_t$ with values
in $B_\h$ clearly has as derivative $iH^\h\times_\h u_t^\h$. We now
adapt to our situation a device which is commonly used to compare
semigroups of operators. As an example quite close to our present
situation, see the proof of equation 16 of
\cite{E}. (Undoubtedly the full expansion of equation 16 could be obtained
in our framework too.)

Fix $t$, and define $\phi$ for this $t$ by $$
\phi(s) =u^\h_s \times_\h (\beta_{t-s}f)\times_\h u_{-s}^\h \ .
$$ From the comments above, $\phi$ is a differentiable function with
values in $B_\h$, whose derivative is given by $$
\phi'(s) = u^\h_s \times_\h ((\pi/i\h)
[H,\beta_{t-s}f]_\h - \{ H,\beta_{t-s}f\}) \times_\h u_{-s}^\h \ .  $$
Notice that $\phi(0)=\beta_tf$ while $\phi(t)=\beta^\h_tf$.  Thus $$
\align
\| \beta^\h_tf-\beta_tf\|_\h  & = \| \int^t_0 \phi'(s)ds\|_\h \\
& \leq |t|\sup \{\| (\pi/i\h)[H,\beta_{t-s}f]_\h -\{ H,\beta_{t-s}f\}
\|_\h: |s| \leq |t|\} \ .
\endalign
$$ From 2.2 above it is clear that it now suffices to control the size
of $$
\| (\pi/ i\h)[H,g]_\h -\{ H,g\}\|_k
$$ where $g=\beta_{t-s}f$.

We now need to use the same arguments as in the proof of theorem 9.3
of \cite{Rf}, but keeping track of $\beta_{t-s}$ so as to get an
estimate which is uniform in $s$. For any multi-index $\mu$ let
$\p^\mu$ denote the corresponding (higher) partial derivative for the
basis $\{E_j\}$ chosen earlier. The norm $\| \ \ \|_k$ is equivalent
to a finite linear combination of the semi-norms $f\mapsto\|\p^\mu
f\|_\infty$ for various $\mu$'s.  So it suffices to obtain suitable
estimates for these semi-norms.  But, just as in the proof of theorem
9.3 of \cite{Rf}, repeated application of Leibniz' rule shows that $$
\|\p^\mu((\pi/i\h)[H,g]_\h -\{ H,g\})\|_\infty
$$ is dominated by a finite linear combination of terms of form $$
\|(\pi/i\h)[\p^\nu H,\p^\lambda g]_\h -
\{ \p^\nu H,\p^\lambda g\}\|_\infty \ ,
$$ where the coefficients of the linear combination do not depend on
$H$, $g$ or $\h$. But $H$ is fixed throughout, and so for notational
simplicity we can set $F=\p^\nu H$ for any given $\nu$. Then we see
that it suffices to obtain for any given multi-index $\lambda$, a
suitable estimate for the size of $$ (\pi/i\h)[F,\p^\lambda g]_\h -\{
F,\p^\lambda g\} \ , $$ where we remember that $g=\beta_{t-s}f$.

To bring all this even closer to the proof of theorem 9.3 of
\cite{Rf}, we use the commutativity of $B$ to write $$ [F,\p^\lambda
g]_\h = (F\times_\h (\p^\lambda g)-F(\p^\lambda g))- ((\p^\lambda
g)\times_\h F-(\p^\lambda g)F) \ .  $$ Then we see that it suffices to
obtain a suitable estimate for the size of $$
(2\pi/i\h)(F\times_\h(\p^\lambda g)- F(\p^\lambda g))-\{ F,\p^\lambda
g\} \tag 2.3 $$ and a similar term.  But by the last displayed
equation in the proof of theorem 9.3 of
\cite{Rf} we find that (2.3) is equal to $\h 2\pi iR(\h)$ where
(after omitting an erroneous subscript $J$) $$ R(\h)=(2\pi i)^{-2}
\sum J_{pj}J_{qk} \int^1_0 \int^1_0 (a_{kj}\times_{\h rs} b_{qp})ds \
rdr \ , $$ and where $a_{kj}=\a_{E_k}\a_{E_j}(F)$ and
$b_{qp}=\a_{E_q}\a_{E_p}(\p^\lambda q)$. From Proposition 2.2 of
\cite{Rf} it then follows that each summand of $R(\h)$ is dominated
for the norm $\| \ \ \|_k$ by $c\|F\|_m\|\p^\lambda g\|_m$ for an
integer $m$ and constant $c$ independent of $F$ and $g$.  Note further
that $\|\p^\lambda g\|_m\leq d\|g\|_n$ for a suitable integer $n$ and
constant $d$. But $\| g\|_n=\|\beta_{t-s}f\|_n$, which is uniformly
bounded for $t-s$ ranging in any finite interval, because of the
continuity given by Theorem 1.4. It follows that $\|R(\h)\|_n$ is
uniformly bounded for $\h\in I$, for our fixed $f$, and for $t$ in any
fixed finite interval.  Because our error term involved $\h R(\h)$, we
thus obtain the desired convergence as $\h\to 0$.  \qed

\bigskip
We remark that one can follow the above analysis more carefully to
obtain a specific bound for $\|R(\h)\|_n$.

We also remark that with somewhat more care we could use the
commutativity of $B$ and the symmetry of $[H,\beta_{t-s}f]$ and $\{
H,\beta_{t-s}f\}$ to obtain an error term of form $\h^2R(\h)$ rather
than $\h R(\h)$, as is usually obtained in discussions of related
situations in the literature, such as expansion 16 of \cite{E}. This
possibility was not discussed in \cite{Rf} since it is not available
when $A$ is not commutative.

The following comments were stimulated by conversations with A.
Vershik. Consider the ordinary 2-torus $T^2$, and let ${\Cal L}_0$
denote $C^\infty(T^2)$ as Lie algebras with the standard Poisson
bracket. The process of associating to elements of ${\Cal L}_0$ their
Hamiltonian vector fields is a Lie algebra homomorphism of ${\Cal
L}_0$ onto the Lie algebra of those smooth vector fields which
generate area-preserving diffeomorphisms of $T^2$. (This homomorphism
is an isomorphism once one factors by the subspace of constant
functions, the center of ${\Cal L}_0$.) As seen in example 10.2 of
\cite{Rf}, the deformation quantization of the symplectic space $T^2$
for the action of ${\Bbb R}^2$ gives the quantum 2-tori (the rotation
algebras) $A_\theta$ (where $\theta = \h$). Let ${\Cal L}_\theta$
denote $A_\theta$ viewed just as a Lie algebra with its commutator
bracket, forgetting the associative algebra structure. It is remarked
in example 3e of \cite{V} (and in references given there and in
\cite{S}) that ${\Cal L}_\theta$ tends to ${\Cal L}_0$ as $\theta$
goes to 0 (with similar statements for other crossed product
algebras). We can view theorem 9.3 of \cite{Rf}, applied to $T^2$, as
then making this intuition rigorous. In the same way, Theorem 2.1 of
the present paper, applied to $T^2$, goes in the direction of saying
rigorously that the group of inner automorphisms of $A_\theta$ coming
from unitaries in the connected component of the unitary group of
$A_\theta$, tends to the group of area-preserving diffeomorphisms of
$T^2$ as $\theta$ goes to 0. (For information on the structure of this
unitary group of $A_\theta$ see \cite{Rf1}.)

\bigskip
\section{3. The proof of Theorem 1.3}

We remark that with suitable care the steps below can be carried out
with $V$ replaced by a general connected Lie group. For simplicity of
exposition we treat only $V$ here since this is all we need, but see the
appendix for the general case.

\proclaim{3.1 LEMMA} Let $N$ be a locally compact space,
let $M=V\times N$, and let $\a$ be the action of $V$ on $M$ coming
from the translation action of $V$ on itself. Let $\Phi\in
B^\infty(M,L)$, and view $\delta_\Phi$ as a derivation of $A^\infty$
(not $B^\infty$). Then $\delta_\Phi$ is the pregenerator of a
one-parameter action $\beta$ on $A$, with corresponding flow $\beta$
on $M$. Furthermore, $\beta$ carries $A^\infty$ into itself, and the
flow $\beta$ on $M$ carries each leaf $V\times\{ n\}$ into itself.
\endproclaim

{\smc Proof}. Because here the action $\a$ of $V$ on $M$ is free, this
lemma is essentially a special case of theorem 2.4.26 of \cite{B}
(which is closely related to results in \cite{BD}). Its proof is
basically just a matter of applying the usual existence theorem for
flows generated by Lipschitz vector fields on $\Bbb R^d$ to obtain a
global flow on each leaf $V\times\{ n\}$.  Then one applies the
theorem concerning the continuous dependence of such flows on their
initial conditions to show that, as $n$ varies, the corresponding
flows fit together continuously to give a flow on $M$. From this
approach we see that $\beta$ carries each leaf of $M$ into itself.
Theorem 2.4.26 of \cite{B} also gives that $A^\infty$ is a core for
the generator of the action $\beta$, and that for each $f\in A^\infty$
we have $$ (d/dt)|_{t=s}\beta_t(f) =\beta_s(\delta_\Phi(f)) \ .  \tag
3.2 $$ Since in our case $\delta_\Phi$ carries $A^\infty$ into itself,
a simple induction argument shows that $\beta$ carries $A^\infty$ into
itself. \qed

\proclaim{3.3 LEMMA} Let $N,M,\a,\Phi$ and $\beta$ be as in the
previous lemma. Let $\beta$ also denote the corresponding
one-parameter action on $C_b(M)$. Then $\beta$ carries $B=B(M,\a)$
into itself, is strongly continuous on $B$, and carries $B^\infty$
into itself. Let $\delta_\Phi$ be the derivation of $B^\infty$ defined
earlier. Then $\delta_\Phi$ is a pregenerator for $\beta$ acting on
$B$.
\endproclaim

{\smc Proof}. Note that we cannot directly invoke theorem 2.4.26 of
\cite{B} here because, in general, the action $\a$ on the maximal
ideal space of $B$ will not be locally free.  In fact, the most
difficult part of the proof is to show that each $\beta_t$ actually
carries $B$ into itself.

Because of the special form of $M$, we can initially work on each leaf
$V\times\{ n\}$ separately. For simplicity of notation we temporarily
consider our $n$ to be fixed, and omit it from the notation, and thus
work on $V$ itself. But we must be careful to obtain estimates which
are uniform in $n\in N$.

By restriction we view $\Phi$ as an element of $B^\infty(V,L)$.  Thus
it is smooth on $V$ in the usual sense. Now $B=B(V)$ will consist
exactly of the uniformly continuous functions on $V$. Thus to show
that $B$ is carried into itself by $\beta$ it clearly suffices to
obtain an estimate of the form $$
\|\beta_t(x)-\beta_t(y)\| \leq K_t\|x-y\|
$$ for all $x,y\in V$, where $K_t$ is a constant independent of $x$
and $y$.  Fix $x,y\in V$ with $x\neq y$, and let $w=y-x$. Let $g$ be
the $V$-valued function on $\Bbb R^2$ defined by $$
g(t,r)=\beta_t(x+rw) \ .  $$ From the usual facts about solutions of
differential equations, $g$ is smooth since $\Phi$ is. Note that for
fixed $t$ the path $g(t,r)$ goes from $\beta_t(x)$ to $\beta_t(y)$ as
$r$ goes from 0 to 1. We consider the length, ${\Cal L}(t)$, of this
path. We use ideas from the first variational equation for ordinary
differential equations (e.g.  page 190 of \cite{A}).

Let $h=\p g/ \p r$, so that $$ {\Cal L}(t) =\int^1_0 \|h(t,r)\| \ dr \
.  $$ Now, by the fact that partial derivatives commute, we have $$
\align
\p h/ \p t & = (\p/ \p r)(\p g/ \p t)
             = (\p/ \p r)(\Phi\circ g) \\ & = ((D\Phi)\circ g)\circ
(\p g/ \p r) = ((D\Phi)\circ g)\circ h \ ,
\endalign
$$ where $D\Phi$ is the usual total derivative of $\Phi$. Note that
since $\beta_t$ is a diffeomorphism and $w\neq 0$, $h$ never takes
value 0, and so the function $\|h(t,r)\|$ is smooth, as is then $\Cal
L$.  A little calculation then shows that $$ |(d{\Cal L}/dt)(t)| \leq
\int^1_0 \|(\p h/ \p t)(t,r)\| dr \leq \| D\Phi\|_\infty \ {\Cal L}(t)
\ .  $$ Consequently $$ {\Cal L}(t) \leq {\Cal L}(0) \
e^{t\|D\Phi\|_\infty} =
\| x-y\| \ e^{t\|D\Phi\|_\infty} \ .
$$ Since ${\Cal L}(t)$ is the length of some curve from $\beta_t(x)$
to $\beta_t(y)$, it follows that $$
\|\beta_t(x)-\beta_t(y)\| \leq \| x-y\| \ e^{t\|D\Phi\|_\infty} \ .
$$ This is an estimate of the desired type, and so as indicated above,
$\beta$ carries $B$ into itself. Notice that we have used the
hypotheses that $\Phi\in B^\infty(V,L)$ to ensure that
$\|D\Phi\|_{\infty}$ is finite.

We return now to the general case in which $M=V\times N$. As long
as we now interpret $\|D\Phi\|_{\infty}$ as a supremum over all of
$M$, which is still finite since $\Phi\in B^\infty(V,L)$, we see that
the above inequality is uniform over all the leaves. It follows easily
that $\beta$ carries $B$ into itself in this case also.

We must now show that the action $\beta$ on $B$ is strongly
continuous. By multiplying elements of $B^\infty$ by elements of
$A^\infty$ which have value 1 on neighborhoods of various points, we
see that every element of $B^\infty$ agrees locally with an element of
$A^\infty$.  It follows that for any $f\in B^\infty$ and any $m\in M$
we have
$$
(d/dt)(f(\beta_t(m))) = (\delta_\Phi f)(\beta_t(m)) \ ,
$$
since this can be viewed as a local statement. In
particular, the derivative on the left exists. Consequently
$$
f(\beta_t(m)) = f(m) + \int^t_0 (\delta_\Phi f)(\beta_s(m))ds \ ,
\tag 3.4
$$
so that $$
\|\beta_t f-f\|_\infty \leq |t| \ \|\delta_\Phi f\|_\infty \ .
$$ Thus $\beta$ is strongly continuous on $B^\infty$. Since $\beta$ is
isometric and $B^\infty$ is dense in $B$, it follows that $\beta$ is
strongly continuous on $B$.

We must now show that $B^\infty$ is contained in the domain of the
infinitesimal generator of $\beta$, and that on $B^\infty$ this
generator agrees with $\delta_\Phi$. We argue much as in the proof of
Lemma 2.4.3 of \cite{B}. Let $f\in B^\infty$, so that $\delta_\Phi
f\in B^\infty$.  Because we now know that $\beta$ is strongly
continuous on $B$, the integral $
\int^t_0 \beta_s (\delta_\Phi f)ds
$ is well-defined for the supremum norm on $B$. Now evaluation at any
point $m\in M$ is continuous for this norm, and so can be brought
inside the integral. From (3.4) it then follows that
$$
\beta_t f-f = \int^t_0 \beta_s(\delta_\Phi f)ds \ .
$$
{}From this it follows immediately that $$ (d/dt)\mid_{t=0} (\beta_t
f) = \delta_\Phi f $$ for the norm on $B$, so that $f$ is in the
domain of the generator of $\beta$. Furthermore, we see that on
$B^\infty$ this generator agrees with $\delta_\Phi$.

It follows readily that equation (3.2) holds for any $f\in B^\infty$.
{}From this equation and the fact that $\delta_\Phi$ carries $B^\infty$
into itself, it follows by a simple induction argument that $\beta$
carries $B^\infty$ into itself. We can now apply Corollary 3.1.7 of
\cite{BR} to conclude that $B^\infty$ is a core for the generator of
$\beta$, i.e. that this generator is the closure of $\delta_\Phi$.
\qed

\bigskip

{\bf 3.5 \ Conclusion of the proof of Theorem 1.3}.  Let $M$,
$\alpha$, $A$ and $B$ be as in the statement of Theorem 1.3.  Let
$P=V\times M$, and let $\tau$ denote the action of $V$ on $P$ coming
from translation on $V$. Let $\eta$ be the map from $P$ to $M$ defined
by $\eta(x,m)=\alpha_x(m)$. Since $\eta$ is surjective, it gives an
isometric isomorphism, still denoted by $\eta$, of $B(M)$ onto a
subalgebra of $B(P)$.  When convenient we will simply identify $B(M)$
with this subalgebra.  Note that $\eta$ is equivariant for $\alpha$
and $\tau$. Thus $B(M)$ is a $\tau$-invariant subalgebra of $B(P)$.
Define (as suggested to me by Alan Weinstein) an action, $\gamma$, of
$V$ on $P$ by $\gamma_y(x,m)=(x-y, \alpha_y(m))$. Note that
$\eta\circ\gamma_y=\eta$ for any $y\in V$, and that the
$\gamma$-orbits of points in $P$ are exactly the $\eta$-preimages of
points in $M$. Let $\Phi\in B^\infty(M,L)$ be given, and set
$\hat\Phi=\Phi\circ\eta$. It is easily seen that $\hat\Phi\in
B^\infty(P,L)$, and clearly $\hat\Phi\circ\gamma_y=\hat\Phi$ for all
$y\in V$. Let $\hat\beta$ denote the flow for $\hat\Phi$ on $P$, whose
existence is assured by Lemma 3.1, and which carries $B^\infty(P)$
into itself by Lemma 3.3. Fix $m\in M$ and $y\in V$. Then $\gamma_y$
gives a bijection of $V\times \{m\}$ onto $V\times\{\alpha_y(m)\}$,
and under this bijection the restrictions of $\hat\Phi$ agree. By the
uniqueness theorem for ordinary differential equations, the
corresponding flows must agree. But by construction these flows are
just given by $\hat\beta$. Thus $\hat\beta$ commutes with each
$\gamma_y$. It follows that for each $t$ the homeomorphism
$\hat\beta_t$ carries each $\gamma$-orbit
into exactly another
$\gamma$-orbit, and so determines a ``flow'', $\beta$, on $M$. It is
easily seen that $\eta$ is an open map ($\gamma$ is a free and proper
action). From this and the continuity of $\hat\beta$
it follows that $\beta$ is continuous, so that it really
is a flow.

Note that by construction $\eta$ is equivariant for $\beta$ and
$\hat\beta$. Since $\eta$ carries $B(M)$ isometrically into $B(P)$ and
$\hat\beta$ carries $B(P)$ into itself and is strongly continuous on
$B(P)$, it follows that $\beta$ carries $B(M)$ into itself and is
strongly continuous there. (Note that $\eta$ does not carry $A(M)$
into $A(P)$.) For the same reasons, $\beta$ will carry $B^\infty(M)$
into itself. A straight-forward calculation using the equivariance of
$\eta$ for $\tau$ and $\alpha$ shows that for $f\in B^\infty(M)$ we
have
$$
\delta_{\hat\Phi}(f\circ\eta)=(\delta_\Phi f)\circ\eta.
$$
Now for fixed $m\in M$ we have $(\beta_t
f)(m)=(\hat\beta_t(f\circ\eta))(0,m)$, and so
$$
(d/dt)|_{t=0}(\beta_t f)(m)=(d/dt)|_{t=0}(\hat\beta_t(f\circ\eta))(0,m)
=(\delta_\Phi f)(m).
$$
We can now argue as in the last parts of the
proof of Lemma 3.3 to conclude that $B^\infty(M)$ is in the domain of
the generator of $\beta$, that on $B^\infty(M)$ this generator agrees
with $\delta_\Phi$, and that $B^\infty(M)$ is a core for this
generator.
\qed

\bigskip
\section{4. The proof of Theorem 1.4}

Exactly as in the conclusion of the proof of Theorem 1.3, let
$P=V\times M$ with action $\tau$ of $V$, so that $B(M)$ is identified
via $\nu$ with a C${}^*$-subalgebra of $B(P)$, and $\beta$ on $B(M)$
is just the restriction of $\hat\beta$ on $B(P)$. Since the action
$\a$ of $V$ on $B(M)$ is just the restriction of the action $\tau$ on
$B(P)$, the semi-norms defined earlier in terms of $\a$ will just be
the restrictions to $B(M)$ of the corresponding semi-norms for $\tau$
on $B(P)$. Thus we see that it suffices to prove Theorem 1.4 for the
setting of Lemma 3.3. This means that it suffices to prove the theorem
on each leaf $V\times\{ m\}$, as long as we obtain uniform estimates
in $m$. Thus we consider first the case $M=V$ with $\a$ the action
$\tau$ of translation, and we consider $\Phi$ and $\beta$ as being on
$V$. Then $\beta$ can be viewed as a function from $\Bbb R\times V$ to
$V$ which is smooth. We will let $D^k$ denote $k$-th derivative for
variables in $V$. We identify $L$ and $V$ in the usual way. Then for
each fixed $(t,x)\in\Bbb R\times V$ the expression $(D^k\beta)(t,x)$
is a symmetric $k$-linear map from $V$ to $V$. We use the inner
product on $V$ to define the norm of this map. The proof of the
following lemma is of a type familiar in the theory of ODE's.

\proclaim{4.1 LEMMA} For any $k > 0$ and any finite interval $I$ about
0 there is a constant $K$ such that $$
\| (D^k\beta)(t,x)\| \leq K
$$ for all $x\in V$ and $t\in I$.
\endproclaim

{\smc Proof}. We argue by induction on $k$. Let $\p$ denote
derivatives with respect to $t$. Thus $\p\beta =\Phi\circ\beta$, where
here and in the following we work pointwise. Now $D$ commutes with
$\p$, so $$
\p(D\beta) = D(\p\beta) =D(\Phi\circ\beta) =
((D\Phi)\circ\beta)\circ D\beta $$ by the chain rule. Thus $$
(D\beta)(t,x) = (D\beta)(0,x) +\int^t_0
(D\Phi(\beta_s(x)))\circ((D\beta)(s,x))ds \ , $$ and so, since
$\beta_0(y)=y$ for all $y$, $$
\|D\beta(t,x)\| \leq 1+\|D\Phi\|_\infty \int^t_0
\|(D\beta)(s,x)\|ds \ .
$$ By Gronwall's inequality (Lemma 4.1.7 of \cite{A}) we obtain $$
\|D\beta(t,x)\| \leq \exp(t\|D\Phi\|_\infty) \ ,
$$ with the right-hand side independent of $x$. Thus the proof of
Lemma 4.1 is complete for $k=1$.

In the proof for higher $k$ we argue by induction. Thus assume that
there is a constant $K_k$ such that for each $j\leq k-1$ we have
$\|(D^j\beta)(t,x)\| \leq K_k$ for $t\in I$ and $x\in V$. Much as
above, we have $$
\p(D^k\beta) = D^k(\p\beta) = D^k(\Phi\circ\beta) \ .
$$ We now invoke the chain rule for higher derivatives, as given for
example on page 92 of \cite{A}. Suppressing some of the notation given
there, we obtain $$
\p(D^k\beta) = ((D\Phi)\circ\beta)\circ (D^k\beta) \ + \ G_k \ \ ,
$$ where $$ G_k = \sum^k_{m=2} \sum ((D^m\Phi)\circ\beta)\circ
(\{(D^{\ell_m}\beta)\}) \ , $$ with each $\ell_m \leq k-1$. By the
induction hypothesis, $\|(D^{\ell_m}\beta)(t,x)\| \leq K_k$ for each
$t\in I$ and $x\in V$.  By our hypothesis on $\Phi$ each
$\|(D^m\Phi)(x)\|$ is uniformly bounded over $V$. It follows that
there is a constant, $L_k$, such that $$
\|G_k(t,x)\| \leq L_k
$$ for all $t\in I$ and $x\in V$. Now $$ (D^k\beta)(t,x) =
(D^k\beta)(0,x)+\int^t_0 (D^k(\Phi\circ\beta))(s,x)ds \ .  $$ But for
$k\geq 2$ we have $(D^k\beta)(0,x)=0$.  Thus $$
\align
\|(D^k\beta)(t,x)\| & \leq \int^t_0 (\|(D\Phi)(\beta_s(x))\| \
\|(D^k\beta)(s,x)\| + \|G_k(s,x)\|)ds \\
& \leq \|D\Phi\|_\infty \int^t_0 \|(D^k\beta)(s,x)\|ds + cL_k \ ,
\endalign
$$ where $c$ is the length of the interval $I$. Thus again by
Gronwall's inequality we find that $$
\|(D^k\beta)(t,x)\| \leq cL_k \exp(t \|D\Phi\|_\infty) \ .
$$
\qed

\bigskip
{\bf Proof of Theorem 1.4}. We deal first with the case $M=V$.  So for
$f\in B^\infty$ we view $D^kf$ as a function from $V$ into the normed
space of $k$-linear maps from $V$ into the complex numbers. Let $I$ be
a finite interval about 0 in $\Bbb R$. As in the above lemma, we let
$\p$ denote derivatives in $t$. Then, working pointwise, we have for
any $f\in B^\infty$ $$
\p(D^k(f\circ\beta)) = D^k(\p(f\circ\beta))
                     = D^k((\delta_\Phi f)\circ\beta) \ .  $$ Set
$g=\delta_\Phi f$, so that $g\in B^\infty$.  Then if we apply the
chain rule to $D^k(g\circ\beta)$, much as in the proof of Lemma 4.1,
and if we apply the conclusion of Lemma 4.1, we find that there is a
constant, $K$, independent of $g$ and $x\in V$, such that $$
\|(D^k(g\circ\beta))(t,x)\| \leq K\left(
\sum_{j\leq k} \| D^jg\|_\infty\right)
$$ for all $t\in I$. Application of the chain rule to
$D^jg=D^j(\delta_\Phi f)$ shows that there is a constant $L$,
depending only on $\Phi$ and $k$, such that $$
\sum_{j\leq k} \|D^jg\|_\infty \leq
L \sum_{j\leq k+1} \|D^jf\|_\infty \ .  $$ Since $$
(D^k(f\circ\beta))(t,x)-(D^kf)(x) =
\int^t_0 D^k(g\circ\beta)(s,x)ds \ ,
$$ it follows that $$
\|(D^k(f\circ\beta_t \, -f))(x)\| \leq
tKL \sum_{j\leq k+1} \|D^jf\|_\infty $$ for all $t\in I$ and $x\in V$.
But the right-hand side is independent of $x$, so we obtain the
desired strong continuity in this case.

Now consider the case where $M=V\times N$ as in Lemma 3.3.  On each
leaf $V\times\{n\}$ we will have the above inequality, and the
constant $KL$ depends only on $\Phi$ and its derivatives on that leaf.
But by examining a bit more carefully the origin of $KL$ and by using
the fact that $\Phi$ and its derivatives are assumed to be uniformly
bounded over all of $M$, we see that we can find finite $KL$ which
works uniformly over all of $M$. The comments at the beginning of this
section complete the proof of strong continuity. The proof of the
remaining facts in the statement of Theorem 1.4 is then essentially
the same as the proof of the similar facts in Lemma 3.3 \qed

\vskip 1 truein
\centerline{\bf Appendix}
\centerline{Lipschitz Flow}

\bigskip
We sketch here how to use a powerful theorem of Derek Robinson
\cite{Rs1} to prove Theorem 1.3. We carry this out in the more
general setting of an action of an arbitrary (connected) Lie group
$G$, with Lie algebra $\frak g$. Let $\a$ be an action of $G$ on a
locally compact space $M$, with corresponding action $\a$ on
$A=C_\infty(M)$, and on $B$, the largest algebra of bounded continuous
functions on $M$ on which $\a$ is strongly continuous. In the same
way, we define $B(M,\frak g)$, a generalization of our earlier
$B(M,L)$. We have the corresponding spaces of smooth vectors
$A^\infty, \ B^\infty$, and $B^\infty(M,\frak g)$.  Much as earlier,
we view elements of $B^\infty(M,\frak g)$ as ``smooth bounded vector
fields'' on $M$.

Let $\Phi\in B^\infty(M,\frak g)$. Then, much as done earlier, $\Phi$
determines a derivation, $\delta_\Phi$, on $A^\infty$ and $B^\infty$.

\proclaim{A1. Theorem}  The derivation $\delta_\Phi$ is the
pregenerator of a one-parameter group, $\beta$, of automorphisms of
$B$, which carries< $A$ into itself, and so determines a flow, $\beta$,
on $M$.
\endproclaim

{\smc Proof}. Robinson's theorem tells us that to show that
$\delta_\Phi$ is a pregenerator it suffices to verify two conditions.
The first is the usual condition that $\delta_\Phi$ be conservative.
The second and crucial condition is that $\delta_\Phi$ be Lipschitz
for the action $\a$.  We recall here briefly what this means.

As we did earlier for $L$, choose an arbitrary inner product on $\frak
g$. This can be translated around $G$ to define a left-invariant
Riemannian metric on $G$. We denote the corresponding length function
\cite{Rs2} by $|x|$. Let $\gamma$ be an action of $G$ on a Banach
space $U$. Let $U^\infty$ denote the smooth vectors for $\gamma$. For
$u\in U^\infty$ let $Du$ denote the linear map from $\frak g$ to $U$
defined by $(Du)(X)=\gamma_X(u)$.  We use the inner product on $\frak
g$ to define $\|Du\|$, and we set $\|u\|_1=\|u\|+\|Du\|$.

\bigskip
{\bf A2. Definition} \cite{Rs1}. With notation as above, an operator
$T:U^\infty\to U$ is said to be Lipschitz if there are constants
$\delta > 0$ and $K$ such that $$
\| {\op{ad}}_{\gamma_x}(T)u\| \leq K |x|\|u\|_1
$$ for all $u\in U^\infty$ and all $x\in G$ with $|x| < \delta$.
(Here \ ad${}_{\gamma_x}(T)=\gamma_x\circ T-T\circ\gamma_x$.)

\bigskip
We will need the following fact, whose proof is a straightforward
argument using the lengths of curves in $G$.

\proclaim{A3. Proposition} With notation as above,
each operator $\gamma_X$ for $X\in\frak g$ is a Lipschitz operator for
$\gamma$.
\endproclaim

We now check that $\delta_\Phi$ and $\a$ satisfy Robinson's
conditions.  That $\delta_\Phi$ is conservative is seen in the usual
way by considering, for any $f\in B^\infty$, the functional consisting
of evaluation at a point of the maximal ideal space of $B$ at which
$f$ takes its maximal absolute value. We now sketch the verification
that $\delta_\Phi$ is a Lipschitz operator for $\a$. It is clear from
the definition that sums of Lipschitz operators are again Lipschitz
operators. Thus it suffices to show that any operator of the form
$h\a_X$ for $h\in B^\infty$ and $X\in\frak g$ is Lipschitz. But such
an operator is the composition of the operators corresponding to
(multiplication by) $h$ and $\a_X$.  It is now convenient for us to
make:

\bigskip
{\bf A4. Definition}. With notation as earlier, we say that $u\in U$
is a Lipschitz vector for the action $\gamma$ if there are constants
$\delta > 0$ and $K$ such that $$
\|\gamma_x(u)-u\| \leq K|x|
$$ for $|x| < \delta$.

\bigskip
Then a straightforward argument again using the length of curves in
$G$ yields the first part of the following proposition. The second
part then follows easily from the first.

\proclaim{A5. Proposition} With notation as earlier,  any $h\in
B^\infty$ is a Lipschitz vector for $\a$. The operator, $M_h$, of
multiplication on $B^\infty$ by $h$ is a Lipschitz vector for
${\op{Ad}}_\a$ and the operator norm.
\endproclaim

We will say that an operator $T$ on $U^\infty$ is of order 1 if there
is an inequality of the form $$
\|Tu\| \leq K\|u\|_1
$$ for $u\in U^\infty$.  We remark that the operators $\a_X$ for
$X\in\frak g$ are clearly of order 1. By a straightforward argument we
then obtain:

\proclaim{A6. Proposition} Let $\gamma$ be an action of $G$ on a
Banach space $U$. Let $T$ be a Lipschitz operator on $U^\infty$ for
$\gamma$, of order 1. Let $S$ be a bounded operator on $U$ which is an
operator-norm Lipschitz vector for ${\op{Ad}}_\gamma$, and carries
$U^\infty$ into itself. Then $ST$ is a Lipschitz operator on
$U^\infty$ for $\gamma$.
\endproclaim

{}From Propositions A3, A5 and A6 it follows that $\delta_\Phi$ is a
Lipschitz operator for $\a$.  We phrased the above discussion for $B$,
but all holds equally well for $A$. We can thus apply Robinson's
theorem to obtain a one-parameter action on $B$ which carries $A$ into
itself, and so gives a flow on $M$.  \qed

\enddocument